\documentclass[11pt, reqno]{amsart}

\newtheorem{teo}{Theorem}[section]
\newtheorem{defin}{Definition}[section]
\newtheorem{prop}{Proposition}[section]
\newtheorem{cor}{Corollary}[section]
\newtheorem{lemma}{Lemma}[section]
\newtheorem{rem}{Remark}[section]

\numberwithin{equation}{section}

\def\proof{{\it Proof:}\ }
\def\endproof{\hfill $\Box$\par\vskip3mm}

\def\eq#1{(\ref{#1})}
\def\pr#1{\ref{#1}}
\def\neweq#1{\begin{equation}\label{#1}}
\def\endeq{\end{equation}}

\def\eps{\varepsilon}
\def\phi{\varphi}

\def\di{\displaystyle}

\def\ep{\varepsilon}
\def\rr{\mathbb{R}}

\def\gr{{\mathcal G}}
\def\fr{{\mathcal F}}

\def\w1{W^{1,1}(\Omega)}

\def\czero{C_0(\overline{\Omega})}

\begin{document}

\title{Decay estimates for a viscous Hamilton-Jacobi equation with
homogeneous Dirichlet boundary conditions} 
\author{Sa\"\i d Benachour}
\address{Institut Elie Cartan, Universit\'e Henri Poincar\'e --
Nancy~I, BP~239, F--54506 Vandoeuvre-l\`es-Nancy, France} 
\email{Said.Benachour@iecn.u-nancy.fr}
\author{Simona D\u abuleanu-Hapca}
\address{SIMBIOS Centre, University of Abertay Dundee, 
DD1 1HG Dundee, UK}
\email{simona.hapca@abertay.ac.uk}
\author{Philippe Lauren\c cot}
\address{Math\'ematiques pour l'Industrie et la Physique,
Universit\'e Paul Sabatier -- Toulouse~3, 118 route de Narbonne,
F--31062 Toulouse cedex 9, France}
\email{laurenco@mip.ups-tlse.fr}
\begin{abstract}
Global classical solutions to the viscous Hamilton-Jacobi equation
$u_t-\Delta u=a\ |\nabla u|^p$ in $(0,\infty)\times\Omega$ with
homogeneous Dirichlet boundary conditions are shown to converge to
zero in $W^{1,\infty}(\Omega)$ at the same speed as the linear heat
semigroup when $p>1$. For $p=1$, an exponential decay to zero is also
obtained in one space dimension but the rate depends on $a$ and
differs from that of the linear heat equation. Finally, if $p\in
(0,1)$ and $a<0$, finite time extinction occurs for non-negative
solutions. 
\end{abstract}

\maketitle

\centerline{\today}

\pagestyle{plain}

\baselineskip16pt

\section{Introduction and main results}

We investigate the large time behaviour of solutions to the following
initial-boundary value problem
\begin{equation}\label{part4.1.1}
\left\{\begin{array}{lll}
\di u_t-\Delta u=a|\nabla u|^p &\mbox{in}& (0,\infty)\times\Omega,\\
\di u=0&\mbox{on}&(0,\infty)\times\partial\Omega,\\
u(0,\cdot)=u_0&\mbox{in}&\Omega ,
\end{array}\right.
\end{equation}
where $a\in\rr$, $a\neq 0$, $p>0$ and $\Omega$ is a bounded
open subset of $\rr^N$ with $C^3$-smooth boundary
$\partial\Omega$. We first recall that several papers have already
been devoted to the well-posedness of the Cauchy-Dirichlet problem
\eq{part4.1.1} 
\cite{Al, BD2, BMP, CLS}. In particular, when the initial datum
is a bounded Radon measure, $p\in[1, (N+2)/(N+1))$ and $a>0$, Alaa proved the
existence and uniqueness of weak solutions to \eq{part4.1.1}
\cite{Al}. When $a>0$ and $p>2$, the non-existence of global solutions is also
studied in \cite{Al,S}, the latter work providing further information
on the way the solution blows up. Using a different approach, Benachour and
Dabuleanu have obtained in \cite{BD2} several results on the existence,
 uniqueness and regularity of global solutions for non-smooth initial data
(typically, $u_0$ is a bounded Radon measure or belongs to
$L^q(\Omega)$ for some $q\ge 1$). These results depend on the sign of
$a$, the value of the exponent $p>0$ and the integrability and sign of
the initial datum $u_0$. Singular initial data had been
considered previously by Crandall, Lions and Souganidis in
\cite{CLS} when $a<0$ and $p>1$: using some properties of
order-preserving semigroups, a universal bound for non-negative solutions
to \eq{part4.1.1} is established in \cite{CLS} which proves useful to show the
existence and uniqueness of solutions to \eq{part4.1.1} when the
initial datum $u_0$ satisfies: $u_0=\infty$ on a bounded open
subset $D\subset\overline{D}\subset\Omega$ and
$u_0=0$ in $\Omega\setminus\overline{D}$.

The main purpose of this paper is to supplement the above mentioned
results by analysing the long time behaviour of global solutions to
\eq{part4.1.1}. While several results are available for the Cauchy
problem \cite{BLS,BLSS,Gi,SZ} and for the Cauchy-Neumann problem
\cite{BD,D}, this question has only been considered recently in
\cite{SZ} for the Cauchy-Dirichlet problem \eq{part4.1.1}: it is shown
there that, for $p>2$, global solutions converge to zero in
$L^\infty(\Omega)$ as time goes to infinity and this property remains
true for global solutions which are bounded in
$C^1(\overline{\Omega})$ when $p\in (1,2]$ (see also \cite{ARBS} for the 
one-dimensional case). Besides giving alternative
proofs of these results, we shall identify the rate at which this
convergence to zero takes place. More precisely, we show that, for
$p>1$, global classical solutions to \eq{part4.1.1} decay to zero in
$W^{1,\infty}(\Omega)$ at the same (exponential) rate as the solutions
to the linear heat equation with homogeneous Dirichlet boundary
conditions. It is only when $p\in (0,1]$ that the gradient term $a
\vert\nabla u\vert^p$ influences the large time dynamics (see
Theorems~\ref{part4.t.1} and~\ref{t.9} below).

\medskip

Before stating our results, we introduce some notations: for $T\in
(0,\infty]$ we set $Q_T=(0,T)\times\Omega~$ and 
$\Gamma_T=(0,T)\times\partial\Omega~$. We denote by
$~C_0(\overline{\Omega})~$ the space of continuous functions on
$\overline{\Omega}$ vanishing on the boundary $\partial\Omega$ and by
$~C_0^+(\overline{\Omega})~$ the positive cone of
$~C_0(\overline{\Omega})~$. Next, $~C^{1,2}(Q_T)~$ is the space of
functions $u\in C(Q_T)$ which are differentiable with respect to
$t\in(0,T)$ and twice differentiable with respect to $x\in\Omega$ with
derivatives $u_t$, $\left( u_{x_i} \right)_{1\leq i\leq N}$ and
$\left( u_{x_i x_j} \right)_{1\leq i,j\leq N }$ belonging to
$C(Q_T)$. For $q\in [1,\infty]$, $\|~ \|_q$ and
$\|~\|_{\partial\Omega,q}$ denote the norms in $L^q(\Omega)$ and
$L^q(\partial\Omega)$, respectively, and $W^{1,q}(\Omega)$ the Sobolev
space of functions in $L^q(\Omega)$ for which the distributional
derivatives $\left( u_{x_i} \right)_{1\leq i\leq N}$ also belong to
$L^q(\Omega)$. Finally, $\nu$ denotes the outward normal unit vector
field to $\Omega$ and we use the notation $u_\nu(t,x) = \nabla
u(t,x)\cdot \nu(x)$ for $(t,x)\in (0,\infty)\times\partial\Omega$ for
the normal trace of the gradient of $u$ (when it is well-defined). 

\medskip

Throughout this paper, we only consider classical solutions to
\eq{part4.1.1} in the following sense:

\begin{defin}\label{d.1}
Given $u_0\in C_0(\overline{\Omega})~$, $p\in (0,\infty)$,
$a\in\rr\setminus\{0\}$ and $T\in (0,\infty]$, a classical solution
$u$ to \eq{part4.1.1} in $Q_T$ is a function $u\in
C([0,T)\times\overline{\Omega})\cap C^{1,2}(Q_T)~$ with $u(0)=u_0$ and
satisfying \eq{part4.1.1} pointwisely in $Q_T$. Such a solution also
satisfies
\begin{equation}\label{1.3}
u(t)=e^{t\Delta} u_0+a\int\limits_0^t e^{(t-s)\Delta}|\nabla u(s)|^p\,ds,\
t\in [0,T),
\end{equation}
where $(e^{t\Delta})_{t\geq 0}$ denotes the semigroup associated to
the linear heat equation with homogeneous Dirichlet boundary
conditions.
\end{defin}

The results concerning the existence of global classical solutions to
\eq{part4.1.1} may then be summarized as follows:

\begin{prop}\label{prexist}
There exists a unique classical solution $u$ to \eq{part4.1.1} in
$Q_\infty$ in the following cases:
\begin{itemize}
\item[(i)] $p\in (0,2]$, $a\in\rr\setminus\{0\}$ and $u_0\in
C_0(\overline{\Omega})~$,
\item[(ii)] $p>2$, $a<0$ and $u_0\in C_0^+(\overline{\Omega})~$,
\item[(iii)] $p>2$, $a>0$ and $u_0\in C_0^+(\overline{\Omega})\cap
C^1(\overline{\Omega})~$ with
$\|u_0\|_{C^1(\overline{\Omega})}\leq\eps_0$ where $\eps_0$ is the
constant defined in \cite[Proposition 3.1]{S}.
\end{itemize}
\end{prop}

We refer to \cite{BD2} for the assertion (i), to
\cite[Theorem~2.1]{CLS} for (ii) and to \cite{S} for (iii). The size
restriction in (iii) is needed to have a global solution since finite
gradient blow-up occurs for sufficiently large initial data
\cite{S}. We also mention that the previous well-posedness statement
is far from being optimal with respect to the regularity of the
initial data. Indeed, on the one hand, existence and uniqueness of
weak solutions to \eq{part4.1.1} (which are classical solutions for
positive times) are established in \cite{BD2} under much weaker
regularity assumptions on $u_0$, depending on the sign of $a$ and the
value of $p$. On the other hand, it has been shown in \cite{BdL04}
that, for any $a\in\rr\setminus\{0\}$, $p>0$ and $u_0\in
C_0(\overline{\Omega})$, the Cauchy-Dirichlet problem \eq{part4.1.1}
has a unique continuous viscosity solution, this result being valid
for nonhomogeneous continuous Dirichlet boundary conditions as
well. Such a solution satisfies the boundary conditions in the
viscosity sense and need not satisfy them in the classical sense for
$p>2$ \cite[p.~62]{BdL04}, the latter phenomenon being in principle
related to the possible blow-up of $\Vert\nabla u(t)\Vert_\infty$. The
restriction on the size of the initial data in
Proposition~\ref{prexist} thus excludes this behaviour. 

\medskip

We have the following results:

\begin{teo}\label{part4.t.1}
Assume that $a<0$, $p\in(0,1)$ and $u_0\in
C_0^+(\overline{\Omega})~$. Denoting by $u$ the corresponding
classical solution to \eq{part4.1.1}, there exists $T^*>0$ such that
\begin{equation}\label{part4.4.19}
u(t,x)=0\ \ ~\mbox{for each}\ \ ~(t,x)\in(T^*,\infty)\times\Omega.
\end{equation}
This property is called \textsl{extinction in finite time} of the
solution to \eq{part4.1.1}.
\end{teo}

The proof relies on the results of \cite{BLS,BLSS,Gi} on the
long time behaviour of the solution to the Cauchy problem in the
whole space  $\rr^N$. Indeed, when $a<0$, after an extension by
$0$ on $\rr^N$ of the initial datum, the solution to the Cauchy
problem becomes a super-solution to the Cauchy-Dirichlet problem.
Thus, from the extinction in finite time of the solutions to the
Cauchy problem we deduce the extinction in finite time of the
solutions to the Cauchy-Dirichlet problem.

\begin{rem}
Extinction in finite time cannot take place for solutions to
\eq{part4.1.1} when $a>0$, $p\in(0,1)$ and $u_0\in
C_0^+(\overline{\Omega})~$, $u_0\not\equiv 0$. Indeed, in this case,
$u$ is greater than the solution to the linear heat equation with
initial datum $u_0$, and thus never vanishes in $\Omega$. Moreover,
when $a>0$ and $p\in(0,1)$, there are non-zero stationary solutions to
\eq{part4.1.1} (see Remark~\pr{part4.r.3}).
\end{rem}

\medskip

In the next result, we establish the convergence to zero of global
solutions to \eq{part4.1.1} for $p\in (1,2]$ and show that it takes
place at the same exponential rate as that of the linear heat
equation.

\begin{teo}\label{part4.t.2}
Assume that $p\in (1,2]$, $a\in\rr\setminus\{0\}$ and $u_0\in
C_0(\overline{\Omega})~$. Denoting by $u$ the corresponding classical
solution to \eq{part4.1.1}, there is a constant $K>0$ depending only
on the initial datum $u_0$, the domain $\Omega$ and the parameters $a$ and
$p$ such that, for $t>0$,
\begin{eqnarray}
\|u(t)\|_{\infty}\leq K e^{-t\lambda_1},\ \ \ \ \qquad\label{4.4}\\
\|\nabla u(t)\|_{\infty}\leq
K(1+t^{-1/2})e^{-t\lambda_1},\label{4.5}
\end{eqnarray}
where $\lambda_1>0$ denotes the first eigenvalue of the Laplace
operator with homogeneous Dirichlet boundary conditions in $\Omega$.
\end{teo}

Recall that, by \cite[Lemma~3, p.~25]{Ro} and \cite[p.~35]{T}, the solution
$e^{t\Delta}u_0$ to the linear heat equation with homogeneous
Dirichlet boundary conditions satisfies:
\begin{eqnarray}
\|e^{t\Delta}u_0\|_{\infty}\leq C_0 e^{-t\lambda_1}\|u_0\|_{\infty},\ \
\qquad\qquad\label{4.6}\\
\|\nabla e^{t\Delta}u_0\|_{\infty}\leq
C_0 (1+t^{-1/2})e^{-t\lambda_1}\|u_0\|_{\infty},\label{4.7}
\end{eqnarray}
for $t>0$ and $u_0\in L^{\infty}(\Omega)$, the constant $C$ depending
only on the domain $\Omega$.
Basically, the previous result asserts that, for $p\in(1,2]$, the
additional nonlinear term, which depends on the gradient, 
has no contribution to the large time behaviour of the solution
whatever the sign of $u_0$ is. Moreover, under the additional
assumptions $a>0$ and $u_0\ge 0$, the temporal decay rate \eq{4.4} is
optimal. Indeed, in that case, the comparison principle ensures that
$u(t)\ge e^{t\Delta}u_0$ and $t\longmapsto \Vert
e^{t\Delta}u_0\Vert_\infty$ behaves as $C\ e^{-t\lambda_1}$ for large
times. More surprisingly, the temporal decay rate \eq{4.4} is also
optimal for a large class of non-negative initial data when $a<0$, so
that the gradient term does not speed up the convergence towards zero
in that case, see Proposition~\ref{prop.1} below. 
As we shall see, a similar
remark is valid when $p>2$ for global solutions to
\eq{part4.1.1} which are bounded in $C^1(\overline{\Omega})$. To prove
Theorem~\ref{part4.t.2}, we will combine the previous decay estimates
for the heat semigroup and a fixed point procedure in an appropriate
weighted Banach space of Kato-Fujita type.

\begin{rem}
Let us point out here that the proof of Theorem~\ref{part4.t.2} for 
$p=2$ is obvious thanks to the Cole-Hopf transformation which reduces 
\eq{part4.1.1} to the linear heat equation (see \eq{spirou} below).
\end{rem}

\medskip

We next show that, if $p>2$, non-negative global solutions to
\eq{part4.1.1} are such that $t\longmapsto e^{t\lambda_1}\Vert
u(t)\Vert_\infty$ is bounded from below and above by positive
constants for large times. More precisely, we have the following
result:

\begin{teo}\label{t.7}
Let $u_0$ be a non-negative function in $C_0^+(\overline{\Omega})\cap
C^1(\overline{\Omega})$ and consider $p>2$. Assume further that either
$a<0$ or $a>0$ and $\|u_0\|_{C^1(\overline{\Omega})}\leq\eps_0$
where $\eps_0$ is the
constant defined in \cite[Proposition 3.1]{S}. Denoting by $u$ the
unique classical solution to \eq{part4.1.1}, there exist non-negative
functions $w_0\in C_0^+(\overline{\Omega})\cap C^1(\overline{\Omega})$
and $W_0\in C_0^+(\overline{\Omega})\cap C^1(\overline{\Omega})$
($w_0\not\equiv 0$ and $W_0\not\equiv 0$ if $u_0\not\equiv 0$)
such that
\begin{equation}\label{4.9}
e^{t\Delta}w_0\leq u(t)\leq e^{t\Delta}W_0 \ \ \ \mbox{ for }\ \
\ t\ \ \ \mbox{ large enough. }
\end{equation}
Moreover, $w_0=u_0$ if $a>0$, $W_0=u_0$ if $a<0$ and $\Vert\nabla
u(t)\Vert_\infty$ satisfies \eq{4.5} .
\end{teo}

\medskip

It turns out that, under the assumptions of either Theorem~\ref{part4.t.2}~(i) 
or Theorem~\ref{t.7}, the estimates \eq{4.4} and \eq{4.5} allow us to be identify 
more precisely the large time behaviour of the solution $u$ to \eq{part4.1.1}. Indeed, 
in these cases, it follows from \eq{4.5} by classical arguments that
$$
\lim_{t\to\infty} \left\Vert e^{\lambda_1 t}\ u(t) - \alpha_\infty\ e_1 \right\Vert_\infty = 0\,, 
$$
where $e_1$ denotes the first eigenfunction of the Laplace operator with homogeneous Dirichlet 
boundary conditions (associated to the eigenvalue $\lambda_1$ and chosen to be non-negative with 
$\Vert e_1\Vert_2=1$) and  
$$
\alpha_\infty = \int_\Omega u_0(x)\ e_1(x)\ dx + a\ \int_0^\infty e^{\lambda_1 t}\ \int_\Omega 
\vert\nabla u(t,x)\vert^p\ e_1(x)\ dxdt\,.
$$
Observe that $\alpha_\infty$ is finite by \eq{4.5} since $p>1$ (but we might have $\alpha_\infty=0$).

\medskip

We finally turn to the case $p=1$ which appears to be a limit
case. Indeed, the proof of the decay rates \eq{4.4} and \eq{4.5}
obtained in Theorem~\ref{part4.t.2} for $p\in (1,2]$ does not extend
to $p=1$. In fact, as we shall see below in the one-dimensional 
case, the $L^\infty$-norm of $u(t)$ still decays exponentially but the
decay rate depends on $a$. The case of several dimensions ($N\geq 2$)
seems to be an open problem. Nevertheless, it has been shown recently in 
\cite{HNR} that there are $\beta\in (0,\lambda_1)$ and 
$U\in C_0^+(\overline{\Omega}) \cap C^2(\overline{\Omega})$ 
such that $-\Delta U = a\ \vert\nabla U\vert + \beta\ U$ in $\Omega$, so that 
$(t,x)\longmapsto e^{-\beta t} U(x)$ is a solution to \eq{part4.1.1}. We
also recall that the case $p=1$ is rather peculiar for the Cauchy
problem in $\rr^N$ \cite{BGL,BRVa,BRVb,GGK,LS}.

\begin{teo}\label{t.9}
Let $\Omega=(-1,1)$, $u_0\in C_0^+(\overline{\Omega})$, $a\in\rr$,
$a\neq 0$ and denote by $u$ the unique classical solution to
\eq{part4.1.1} with $p=1$. Then there is a positive constant $\gamma(a)$
depending only on $a$ and $\Omega$ such that $u$ satisfies
\begin{equation}\label{7.2}
\|u(t)\|_{\infty}\leq \gamma(a)\ \|u_0\|_{\infty}\
e^{-((a^2/4)+\alpha_1)t}, \quad t\ge 1,
\end{equation}
where $\alpha_1$ is the first eigenvalue of the unbounded linear
operator $L$ of $L^2(0,1)$
with domain
$$
D(L) = \left\{ \varphi\in W^{2,2}(0,1) \;\;\mbox{ such that }\;\;
a\varphi(0)+2\varphi_x(0)=\varphi(1)=0 \right\}
$$
and defined by
$$
L(\varphi)=-\varphi_{xx} \;\;\mbox{ for }\;\; \varphi\in D(L).
$$
\end{teo}

The restriction to the one-dimensional setting comes from the
observation that, if $u_0$ is an even function in $(-1,1)$ which is
non-increasing in $(0,1)$, then $u$ also solves
$$
\left\{\begin{array}{lll} u_t-u_{xx}=-a
u_x&\mbox{in}&(0,\infty)\times(0,1),\\
u_x(t,0)=u(t,1)=0&\mbox{in}&(0,\infty),\\
u(0,x)=u_0(x)&\mbox{in}&(0,1),\end{array}\right.
$$
which is a linear convection-diffusion equation. After a suitable
change of unknown function, the study of the large time behaviour of
$u$ reduces to the spectral decomposition in $L^2(0,1)$
of $L$ whose eigenfunctions form an orthonormal basis in
$L^2(0,1)$. The case of general initial data will then follow by a
comparison argument.

\begin{rem}
It follows from the previous analysis that the gradient term $\vert\nabla u\vert^p$ 
alters the large time dynamics only for $p=1$ which contrasts markedly with the Cauchy 
problem in $\rr^N$ where the effects of the gradient term become preponderant for 
$p<(N+2)/(N+1)$, see \cite{BLS,BLSS,Gi,LS} and the references therein.   
\end{rem}

%%%%%%%%%%%%%%%%%%%%%%%%%%%%%%%%%%%%%%%%%%%%%%%%%%%%%%%%%%%%%%%%%%%%%%%

\section{Proof of Theorem \pr{part4.t.1}}

Let $\tilde{u}_0$ be the extension by $0$ of $u_0$ outside
the domain $\Omega$, that is, $\tilde{u}_0(x)=u_0(x)$ if $x\in\Omega$
and $\tilde{u}_0=0$ if $x\in\rr^N\setminus\Omega$, and denote by
$\tilde{u}\in C^{1,2}((0,\infty)\times\rr^{N})\cap
C([0,\infty)\times\rr^N)$ the unique solution to the Cauchy problem
\cite{GGK}
\begin{equation}\label{part4.4.26}
\left\{\begin{array}{ll}
\di\tilde{u}_t-\Delta\tilde{u}=a|\nabla\tilde{u}|^p\ \ ~\mbox{in}\ \
~(0,\infty)\times\rr^N,\\
\tilde{u}(0,.)=\tilde{u}_0\ \ ~\mbox{in}\ \ ~\rr^N.\end{array}\right.
\end{equation}
Since $a<0$,  $p\in(0,1)$ and $\tilde{u}_0$ is a non-negative
continuous and bounded function with compact support in $\rr^N$, we
infer from \cite{BLS,Gi} that $\tilde{u}$ enjoys the property of
extinction in finite time, that is, there exists $T^*>0$ such that
$\tilde{u}(t,x)=0$ for $(t,x)\in(T^*,\infty)\times\rr^N$. On the
other hand, $\tilde{u}$ is non-negative by the comparison principle
and thus satisfies:
\begin{equation}\label{part4.4.17}
\left\{\begin{array}{lll}
\di\tilde{u}_t-\Delta\tilde{u}=a|\nabla\tilde{u}|^p\ \ ~\mbox{in}\ \
~Q_\infty,\\
\tilde{u}(t,x)\geq 0\ \ \ \ \mbox{on}\ \ ~\Gamma_{\infty},\\
\tilde{u}(0,.)=u_0\ \ ~\mbox{in}\ \ ~\Omega.
\end{array}\right.\end{equation}
Therefore, $\tilde{u}$ is a super-solution to \eq{part4.1.1} and the
comparison principle \cite[Theorem~16, p.~52]{F} ensures that $0\leq
u\leq\tilde{u}$ in $Q_{\infty}$. Consequently, there is $T^*>0$ such
that $u$ satisfies \eq{part4.4.19}.
\endproof

%%%%%%%%%%%%%%%%%%%%%%%%%%%%%%%%%%%%%%%%%%%%%%%%%%%%%%%%%%%%%%%%%%%%%%%%

\section{Convergence to zero for $p\in (1,2]$}

As a preliminary step to the proof of Theorem~\ref{part4.t.2} we first
establish the convergence to zero in $L^\infty(\Omega)$ of any global
classical solution to \eq{part4.1.1} when $p\in (1,2]$. 

\begin{prop}\label{part4.p.1}
Assume that $p\in (1,2]$, $a\in\rr\setminus\{0\}$ and $u_0\in
C_0(\overline{\Omega})$. Denoting by $u$ the corresponding classical
solution to \eq{part4.1.1} we have
$$
\lim_{t\to\infty} \Vert u(t)\Vert_\infty=0\,.
$$
\end{prop}

A possible proof of Proposition~\ref{part4.p.1} relies on the
LaSalle Invariance Principle since it can be shown that the
$L^\infty$-norm is a strict Liapunov functional for the dynamical
system associated to \eq{part4.1.1} in
$C_0(\overline{\Omega})$. However, the following shorter proof relying
on the method of relaxed semi-limits in the spirit of \cite[Section~8,
Exemple~5]{Ba97} has been suggested to us by G.~Barles
\cite{Ba06}.  

\proof\\ We first introduce some notations: for $\ep\in [0,1)$,
$\xi_0\in\rr$, $\xi=(\xi_i)_{\{1\le i\le N\}}\in\rr^N$ and any
symmetric $N\times N$ matrix $S\in\mathcal{M}_N(\rr)$, we put 
$$
G^\ep(t,x,r,\xi_0,\xi,S) = \left\{ 
\begin{array}{lcl}
\ep\ \xi_0 - tr(S) - a \ \vert\xi\vert^p & \mbox{ if } & (t,x)\in
\mathcal{O}\,, \\
r & \mbox{ if } & (t,x)\in (0,\infty)\times\partial\Omega\,,\\
r - u_0(x) & \mbox{ if } & (t,x)\in \{0\}\times\overline{\Omega}\,,
\end{array}
\right.
$$
and
$$
H(x,r,\xi,S) = \left\{ 
\begin{array}{lcl}
- tr(S) - a \ \vert\xi\vert^p & \mbox{ if } & x\in \Omega\,, \\
r & \mbox{ if } & x\in \partial\Omega\,,
\end{array}
\right.
$$ 
where $tr(S)$ denotes the trace of the matrix $S$ and
$\mathcal{O}=(0,\infty)\times\Omega$.

\medskip  

Next, we assume that $a<0$ and put $u^\ep(t,x) = u(t/\ep,x)$ for
$\ep\in (0,1)$ and $(t,x)\in \overline{\mathcal{O}}$. It
readily follows from \eq{part4.1.1} and Proposition~\ref{prexist} that
$u^\ep\in C(\overline{\mathcal{O}})$ and solves 
\begin{equation}
\label{bar1}
G^\ep(t,x,u^\ep,u_t^\ep,\nabla u^\ep,D^2 u^\ep) = 0 \;\;\mbox{ in
}\;\; \overline{\mathcal{O}}
\end{equation}
in the viscosity sense ($D^2 u^\ep$ denoting the Hessian matrix of
$u^\ep$). In addition, the following bound is a
straightforward consequence of the maximum principle
\begin{equation}
\label{bar2}
\Vert u^\ep(t) \Vert_\infty \le \Vert u_0\Vert_\infty\,, \quad
(t,\ep)\in [0,\infty)\times (0,1)\,.
\end{equation}
Introducing the semi-limits $\underline{u}$ and $\overline{u}$ defined
by
$$
\underline{u}(t,x) = \liminf_{(\ep,s,y)\to (0,t,x)} u^\ep(s,y)\,,
\quad \overline{u}(t,x) = \limsup_{(\ep,s,y)\to (0,t,x)} u^\ep(s,y)\,,
$$
for $(t,x)\in \overline{\mathcal{O}}$, we clearly have
\begin{equation}
\label{bar3}
\underline{u}(t,x) \le \overline{u}(t,x)\,, \quad
(t,x)\in\overline{\mathcal{O}}\,, 
\end{equation}
and infer from (\ref{bar1}), (\ref{bar2}) and \cite[Th\'eor\`eme~
4.1]{Ba94} that $\overline{u}$ is an upper semicontinuous viscosity
subsolution to $G^0=0$ in $\overline{\mathcal{O}}$ while
$\underline{u}$ is a lower semicontinuous viscosity supersolution to
$G^0=0$ in $\overline{\mathcal{O}}$. Since $u^\ep$ is obtained from
$u$ by a time dilatation, we realize that $\underline{u}$ and
$\overline{u}$ actually do not depend on the time variable, that is,  
\begin{equation}
\label{bar4}
\underline{u}(t,x) = \underline{u}(1,x) = \underline{v}(x) \;\;\mbox{
and }\;\; \overline{u}(t,x) = \overline{u}(1,x) = \overline{v}(x)
\;\;\mbox{ for }\;\; (t,x)\in (0,\infty)\times\overline{\Omega}\,.
\end{equation}
We therefore deduce from the properties of $\underline{u}$ and
$\overline{u}$ that $\overline{v}$ is an upper semicontinuous
viscosity subsolution to $H=0$ in $\overline{\Omega}$ while
$\underline{v}$ is a lower semicontinuous viscosity supersolution to
$H=0$ in $\overline{\Omega}$. In addition, since the boundary of
$\Omega$ is smooth, we may proceed as in the proof of
\cite[Corollary~2.1]{BaBu01} or apply \cite[Theorem~2.1]{Gri06} to
deduce that 
\begin{equation}
\label{bar5}
\overline{v}(x) \le 0 \le \underline{v}(x) \;\;\mbox{ for }\;\;
x\in\partial\Omega\,. 
\end{equation}
Indeed, $H$ clearly fulfils \cite[(F5)]{BaBu01} and
\cite[Theorem~2.1~(iv)]{Gri06}. We are now in a position to apply the
strong comparison principle stated in Proposition~\ref{scr} below to
conclude that 
\begin{equation}
\label{bar6}
\overline{v}(x) \le \underline{v}(x) \;\;\mbox{ for }\;\;
x\in\overline{\Omega}\,. 
\end{equation} 
Combining \eq{bar3}, \eq{bar4} and \eq{bar6}, we conclude that
$\overline{v}=\underline{v}=v$ in $\overline{\Omega}$ and, thanks to
\cite[Lemma~4.1]{Ba94}, we obtain that $v\in C(\overline{\Omega})$ is
a viscosity solution to $H=0$ and $(u^\ep)$ converges towards $v$ in
$C((0,\infty)\times\overline{\Omega})$. In particular, $(u^\ep(1,.))$
converges towards $v$ in $C(\overline{\Omega})$ from which we deduce
that
$$
\lim_{t\to\infty} \Vert u(t) - v \Vert_\infty = 0\,.
$$
On the other hand, $x\mapsto 0$ is also a continuous viscosity
solution to $H=0$ and it readily follows from \eq{bar5} and
Proposition~\ref{scr} that $v=0$, which completes the proof of
Proposition~\ref{part4.p.1} when $a<0$. 

Finally, if $a>0$, it is straightforward to check that $-u$ solves
\eq{part4.1.1} with $-a$ instead of $a$ and $-u_0$ instead of
$u_0$. We then apply the previous analysis to $-u$ to complete the
proof of Proposition~\ref{part4.p.1}. \endproof

\medskip

We now turn to the cornerstone of the previous proof, namely the
strong comparison principle \cite[Theorem~3.3]{CIL92}. We cannot
however apply directly \cite[Theorem~3.3]{CIL92} because $H$ lacks
some coercivity with respect to $r$ and thus does not fulfil
\cite[(3.13)]{CIL92}. Still, the fact that \cite[Theorem~3.3]{CIL92}
holds true without this assumption in some cases is already well-known
and a strategy to bypass this assumption is sketched in
\cite[Section~5.C]{CIL92} (see also \cite[Section~4.4.1]{Ba97}). We
will thus only give the required modification of the proof of
\cite[Theorem~3.3]{CIL92}. Actually, Proposition~\ref{scr} is a
particular case of \cite[Theorem~1]{BaBu01} which applies to a wider
class of functions $H$ with no dependence on $u$ but is more
complicated to prove.     

\begin{prop}\label{scr}
Let $U$ be an upper semicontinuous viscosity subsolution to $H=0$ in
$\overline{\Omega}$ and $V$ be a lower semicontinuous viscosity
supersolution to $H=0$ in $\overline{\Omega}$ such that $U\le V$ on
$\partial\Omega$. If $a<0$, then $U\le V$ in $\overline{\Omega}$.  
\end{prop}

\proof\\ According to \cite[Section~5.C]{CIL92} and
\cite[Section~4.4.1]{Ba97}, it suffices to construct a sequence
$(U_\delta)_{\delta\in (0,1)}$ such that (i) $U_\delta\le V$ on
$\partial\Omega$, (ii) $U_\delta$ is an upper semicontinuous viscosity
subsolution to $H+\eta(\delta)=0$ for some $\eta(\delta)>0$, (iii)
$\eta(\delta)\to 0$ and $\Vert U_\delta- U\Vert_\infty\to 0$ as
$\delta\to 0$. 

We now construct such a sequence: since $\Omega$ is bounded, there
exists $\lambda>0$ such that $e^{(p-1)(x_1-\lambda)} \le 1/(2\vert
a\vert)$. We then put 
$$
M = \max_{x\in\overline{\Omega}} e^{x_1-\lambda}\,, \quad m =
\min_{x\in\overline{\Omega}} e^{x_1-\lambda} 
$$
and
$$
U_\delta(x) = (1-\delta)\ U(x) + \delta\ \left( e^{x_1-\lambda} - M +
\min_{y\in\overline{\Omega}} V(y) \right)\,, \quad
x\in\overline{\Omega} \,. 
$$
Owing to the convexity of $H$ with respect to $\xi$ and $S$ (recall
that $a<0$), $U_\delta$ satisfies the requirements (i), (ii) and (iii)
listed above with $\eta(\delta)=(m\ \delta)/2$. At this point, one
then proceeds as in the proof of \cite[Section~5.C]{CIL92} to conclude
that $U_\delta\le V$ in $\overline{\Omega}$ for each $\delta\in (0,1)$
and then pass to the limit as $\delta\to 0$ to complete the proof. 
\endproof

\medskip

\begin{rem}\label{part4.r.3}
It follows from Proposition~\ref{part4.p.1} that zero is the only
bounded stationary solution to \eq{part4.1.1} when $p\in (1,2]$. This
property no longer holds true for
$p\in(0,1)$. Indeed, in the particular case where $\Omega$ is the unit
ball $B_1(0)$ of $\rr^N$ and $a=1$, there are at least two stationary
solutions to \eq{part4.1.1} as the function
$$w(x) = \frac{(1-p)^{(2-p)/(1-p)}}{(2-p)(N-(N-1)p)^{1/(1-p)}}
\left( 1-|x|^{(2-p)/(1-p)} \right), \quad x\in B_1(0),$$
is a non-zero stationary solution to \eq{part4.1.1} in
$\Omega=B_1(0)$ (for $p>1$, similar solutions exist but are
singular \cite{AP,Lions}). Furthermore, if $N=1$, there is a continuum
of non-negative stationary solutions and the convergence towards these
stationary solutions is investigated in \cite{La}. For general domains
$\Omega$ in several space dimensions, the large time behaviour seems
to be an open problem.
\end{rem}

%%%%%%%%%%%%%%%%%%%%%%%%%%%%%%%%%%%%%%%%%%%%%%%%%%%%%%%%%%%%%%%%%%%%%%%

\section{Proof of Theorem \pr{part4.t.2} - Decay estimates}
We denote by $X$ the Banach space of functions
$u\in C([0,\infty),C_0(\overline{\Omega}))\cap
C((0,\infty),C^1(\overline{\Omega}))$ 
such that
$$
\Vert u\Vert_X = \max\left\{ \sup\limits_{t\in(0,\infty)}
e^{t\lambda_1} \|u(t)\|_{\infty} , \sup\limits_{t\in(0,\infty)}
\frac{t^{1/2}}{1+t^{1/2}} e^{t\lambda_1}\|\nabla
u(t)\|_{\infty}\right\}<\infty,
$$
where $\lambda_1$ is the first eigenvalue of the Laplace operator with
homogeneous Dirichlet boundary conditions in $\Omega$. Similar
weighted spaces related to the heat semigroup have been previously
used by, e.g., Kato \cite{Kt}, Kato and Fujita \cite{KF}, Brezis and Cazenave
\cite{BC} and Ben-Artzi, Souplet and Weissler \cite{BSW}.

For $u\in X$, $u_0\in\czero$ and $t>0$, we introduce the maps
\begin{equation}\label{3.4.1}
{\mathcal{G}}u(t)=\int\limits_{0}^{t} e^{(t-s)\Delta}|\nabla
u(s)|^p\,ds
\end{equation}
 and
\begin{equation}\label{3.4.2}
{\mathcal{F}}u(t)=e^{t\Delta}u_0+a {\mathcal{G}}u(t).
\end{equation}

Clearly, a solution to \eq{part4.1.1} is a fixed point of
$\mathcal{F}$ and we will use a fixed point procedure to show that
some solutions to \eq{part4.1.1} belong to some suitable bounded
subsets of $X$. More precisely, we have the following result:

\begin{prop}\label{fantasio}
If $p\in (1,2)$, there is a positive constant $K_1>0$ depending
only on $p$, $a$,
$\lambda_1$ and $\Omega$ such that, if $u_0\in \czero$ satisfies
$\Vert u_0\Vert_\infty\le K/(2C_0)$ for some $K\in (0,K_1]$, the
corresponding classical solution $u$ to \eq{part4.1.1} belongs to
$X_K=\{ u\in X; \Vert u\Vert_X\le K\}$ (recall that $C_0$ is the
constant occurring in \eq{4.6} and \eq{4.7}).
\end{prop}

\proof\\ Recalling that the Laplace operator with homogeneous Dirichlet
boundary conditions in $\Omega$ generates an analytic semigroup in
$\czero$ and $C_0(\overline{\Omega})\cap C^1(\overline{\Omega})$ (see,
e.g., \cite[Definition~2.0.2, Corollary~3.1.21 \& Theorem~3.1.25]{Lu}),
both $\mathcal{G}u$ and $\mathcal{F}u$ belong to
$C([0,\infty),C_0(\overline{\Omega}))\cap
C((0,\infty),C^1(\overline{\Omega}))$ for
$u\in X$. We next check that $\mathcal{G}u$ and $\mathcal{F}u$ map $X$
into itself. Let $u\in X$. On the one hand, taking into account
\eq{4.6} we have
$$\begin{array}{l}
\|\gr
u(t)\|_{\infty}\leq\di\int\limits_0^t\|e^{(t-s)\Delta}|\nabla
u(s)|^p\|_{\infty}\,ds\leq
C_0\di\int\limits_0^te^{-(t-s)\lambda_1}\|\nabla
u(s)\|^p_{\infty}\,ds\\
\ \ \ \ \ \ \ \ \ \ \ \ \leq C_0 \|u\|^p_{X}\di\int\limits_0^t
e^{-(t-s)\lambda_1} (1+s^{-1/2})^p e^{-sp\lambda_1}\,ds\\
\ \ \ \ \ \ \ \ \ \ \ \ \leq C_0 \|u\|^p_{X} e^{-t\lambda_1}
\di\int\limits_0^t (1+s^{-1/2})^p e^{-s(p-1)\lambda_1}\,ds,\\
\end{array}$$
which implies
\begin{equation}\label{6.1}
e^{t\lambda_1}\|\gr u(t)\|_{\infty}\leq C_0\ I_1(p)
\|u\|^p_{X} \;\;\mbox{ with }\;\; I_1(p) = \di\int\limits_0^\infty
(1+s^{-1/2})^p e^{-s(p-1)\lambda_1}\,ds<\infty,
\end{equation}
the integral $I_1(p)$ being finite since $p\in (1,2)$. On the
other hand, we infer from \eq{4.7} that
$$\frac{t^{1/2}}{1+t^{1/2}}e^{t\lambda_1}\|\nabla\gr
u(t)\|_{\infty}\leq C_0
\|u\|^p_{X}\frac{t^{1/2}}{1+t^{1/2}}\di\int\limits_0^t
(1+(t-s)^{-1/2})(1+s^{-1/2})^p e^{-s(p-1)\lambda_1}\,ds.$$
Since $p\in (1,2)$, we have
\begin{eqnarray*}
& & \frac{t^{1/2}}{1+t^{1/2}}\di\int\limits_0^t
(1+(t-s)^{-1/2})(1+s^{-1/2})^p e^{-s(p-1)\lambda_1}\,ds \\
& \le & \frac{2t^{1/2}}{1+t^{1/2}}\di\int\limits_0^t \left( 1 +
(t-s)^{-1/2} + s^{-p/2} + (t-s)^{-1/2}\ s^{-p/2} \right)
e^{-s(p-1)\lambda_1}\,ds \\
& \le & 2 \di\int\limits_0^t \left( 1 + (t-s)^{-1/2} + s^{-p/2} +
\frac{t^{1/2}}{1+t^{1/2}} (t-s)^{-1/2}\ s^{-p/2} \right)
e^{-s(p-1)\lambda_1}\,ds \\
& \le & 2 \di\int\limits_0^\infty \left( 1 +s^{-p/2} \right)
e^{-s(p-1)\lambda_1}\,ds + 2 \sup_{r\ge 0}\left\{ r^{1/2}\
e^{-r(p-1)\lambda_1} \right\} \di\int\limits_0^1 (1-s)^{-1/2}
s^{-1/2}\ ds  \\
& + & 2 \frac{t^{(2-p)/4}}{1+t^{1/2}}\ \sup_{r\ge 0} \left\{
r^{(2-p)/4}\ e^{-r(p-1)\lambda_1} \right\} \di\int\limits_0^1
(1-s)^{-1/2}\ s^{-(p+2)/4}\ ds,
\end{eqnarray*}
and the right-hand side of the above inequality is bounded since $p\in
(1,2)$ and $t^{(2-p)/4}\le 1+t^{1/2}$ for $t\ge 0$. Consequently,
\begin{equation}
\label{spip}
I_2(p) = \sup_{t\ge 0} \left\{
\frac{t^{1/2}}{1+t^{1/2}}\di\int\limits_0^t
(1+(t-s)^{-1/2})(1+s^{-1/2})^p e^{-s(p-1)\lambda_1}\,ds \right\} <
\infty
\end{equation}
and
\begin{equation}\label{6.2}
\frac{t^{1/2}}{1+t^{1/2}}e^{t\lambda_1}\|\nabla\gr
u(t)\|_{\infty}\leq C_0 I_2(p)\|u\|_X^p.
\end{equation}
Combining \eq{6.1} and \eq{6.2} we conclude that
$$\|\gr(u)\|_X\leq C_2\|u\|_X^p,$$
where $C_2$ is a constant depending only on $p$, $\lambda_1$ and
$\Omega$. Applying once more the estimates \eq{4.6} and \eq{4.7} and
using the above inequality we obtain
\begin{equation}\label{6.3}
\|\fr(u)\|_X\leq C_0\ \|u_0\|_{\infty}+\vert a\vert C_2 \|u\|_X^p.
\end{equation}
We have thus established that $\mathcal{F}u\in X$ for $u\in X$.

Next, for $K>0$, we consider two functions $u_1$ and  $u_2$ in
$X_K$, where $X_K$ is defined in Proposition~\ref{fantasio}. 
By \eq{4.6} we have
\begin{eqnarray*}
\|\fr u_1(t)-\fr u_2(t)\|_{\infty} & \leq & \vert
a\vert\di\int\limits_0^t\|e^{(t-s)\Delta}(|\nabla u_1|^p(s)-|\nabla
u_2|^p(s))\|_{\infty}\, ds\\
& \leq & \vert a\vert p C_0 \di\int\limits_0^t
e^{-(t-s)\lambda_1}(\|\nabla u_1(s)\|_{\infty}^{p-1}+ \|\nabla
u_2(s)\|_{\infty}^{p-1})\|\nabla (u_1-u_2)(s)\|_{\infty}\,ds \\
& \le & \vert a\vert p C_0 I_1(p) (\|u_1\|_X^{p-1}+
\|u_2\|_X^{p-1})\Vert u_1 - u_2\Vert_X e^{-t\lambda_1}\,,
\end{eqnarray*}
whence
\begin{equation}\label{6.4} e^{t\lambda_1}\|\fr
u_1(t)-\fr u_2(t)\|_{\infty}\leq C_3 K^{p-1} \|u_1-u_2\|_X,
\end{equation}
where $C_3$ is a constant depending only on $p$, $\lambda_1$ and
$\Omega$. Likewise we deduce from \eq{4.7} and \eq{spip} that
\begin{equation}\label{6.5}
\frac{t^{1/2}}{1+t^{1/2}}e^{t\lambda_1}\|\nabla\fr
u_1(t)-\nabla\fr u_2(t)\|_{\infty}\leq C_3 K^{p-1} \|u_1-u_2\|_X
\end{equation}
for a possibly larger constant $C_3$. Combining \eq{6.4} and \eq{6.5}
the functional $\fr$ satisfies
\begin{equation}\label{6.6} \|\fr u_1-\fr u_2\|_{X}\leq
C_3 K^{p-1}\|u_1-u_2\|_X
\end{equation}
for $u_1\in X_K$ and $u_2 \in X_K$. Introducing $K_1>0$ given by
\begin{equation}
\label{Kcrit}
K_1^{p-1} = \min{\left\{ \frac{1}{2 \vert a\vert C_2} , \frac{1}{2C_3}
\right\}},
\end{equation}
we infer from \eq{6.3} and \eq{6.6} that, if $K\in (0,K_1]$ and $\Vert
u_0\Vert_\infty\le K/(2C_0)$, then $\fr(u)\in X_K$ for $u\in X_K$ and
$\|\fr u_1-\fr u_2\|_X\leq \|u_1-u_2\|_X~/2~$ for $u_1\in X_K$ and
$u_2 \in 
X_K$. Consequently, under these assumptions on $K$ and $u_0$, $\fr$ is
a strict contraction from $X_K$ into $X_K$. By the Banach fixed point
theorem, $\fr$ has a unique fixed point $u\in X_K$. Then $u$ satisfies
\eq{1.3} and \cite[Theorems~3.2 \&~3.3]{BD2}  warrant that $u$ is the
unique classical solution to \eq{part4.1.1}, which completes the proof
of Proposition~\ref{fantasio}.
\endproof

\medskip

{\it Proof of Theorem \pr{part4.t.2} (decay estimates):}\\
Assume first that $p\in (1,2)$. We consider $u_0\in\czero$ and
denote by $u$ the corresponding classical
solution to \eq{part4.1.1}. We have already established in the
previous section that $\Vert u(t)\Vert_\infty \longrightarrow 0$ as
$t\to \infty$. Therefore, there exists $t_0>0$
such that $\|u(t_0)\|_{\infty}\le K_1/(2C_0)$ and we infer from
Proposition~\ref{fantasio} with $K=K_1$ that $u(.+t_0)$ belongs to
$X_{K_1}$, i.e.,
$$
\|u(t)\|_{\infty}\leq K_1 e^{-(t-t_0)\lambda_1} \;\;\mbox{ and }\;\;
\|\nabla u(t)\|_{\infty}\leq K_1
(1+(t-t_0)^{-1/2}e^{-(t-t_0)\lambda_1}
$$
for any $t>t_0$. On the other hand, it follows from the analysis in
\cite{BD2} that $\Vert u(t)\Vert_\infty \le \Vert u_0\Vert_\infty$ and
$\|\nabla u(t)\|_{\infty}\leq C(t_0,u_0)t^{-1/2}$ for $t\in
(0,t_0]$. Combining these two facts yields \eq{4.4} and \eq{4.5} for
$u$.

\medskip

It remains to study the case $p=2$. Introducing $U=e^{a u} -1$, it
follows from \eq{part4.1.1} that $U_t=\Delta U$ in
$(0,\infty)\times\Omega$ with $U=0$ on
$(0,\infty)\times\partial\Omega$ and $U(0)=U_0=e^{a u_0} -1$ in
$\Omega$. Consequently, $U(t)=e^{t\Delta} U_0$ and
\begin{equation}
\label{spirou}
u(t) = \frac{1}{a}\ \log{\left( 1 + e^{t\Delta} U_0 \right)}
\;\;\mbox{ and }\;\; \nabla u(t) = \frac{1}{a}\ \frac{\nabla 
e^{t\Delta} U_0}{1 + e^{t\Delta} U_0}\,, \quad
t\ge 0\,.
\end{equation}
Since $\Vert e^{t\Delta}
U_0\Vert_\infty\longrightarrow 0$ as $t\to\infty$, the temporal
decay estimates \eq{4.4} and \eq{4.5} follow from \eq{4.6} and
\eq{4.7} and the previous formulae for $u(t)$ and $\nabla u(t)$.
\endproof

\medskip

As already mentioned, when $a>0$, $p\in (1,2]$ and $u_0\in
C_0^+(\overline{\Omega})$, the temporal decay rate \eq{4.4} is optimal
since $u(t)\ge e^{t\Delta}u_0$ by the comparison principle and
$t\longmapsto \Vert e^{t\Delta}u_0\Vert_\infty$ behaves as $C\
e^{-t\lambda_1}$ for large times. It turns out that the temporal decay
rate \eq{4.4} is still optimal when $a<0$, $p\in (1,2)$ and $u_0\in
C_0^+(\overline{\Omega})$.

\begin{prop}\label{prop.1}
Consider $a<0$, $p\in(1,2)$ and $u_0\in C_0^+(\overline{\Omega})$ such
that $u_0\geq\alpha e_1$ for some $\alpha>0$, where $e_1$ denotes the
eigenfunction of the Laplace operator with homogeneous Dirichlet
boundary conditions in $\Omega$ associated to the eigenvalue
$\lambda_1$ and normalized such that $e_1\geq 0$ and
$\|e_1\|_{\infty}=1$. Then, there exist two positive constants
$C_4(u_0)>0$ and $C_5(u_0)>0$ depending on $u_0$ such that the
classical solution $u$ to \eq{part4.1.1} satisfies:
\begin{equation}\label{6.10}
C_4(u_0)e^{-t\lambda_1}\leq \|u(t)\|_{\infty}\leq
C_5(u_0)e^{-t\lambda_1}\ \ \ \mbox{for}\ \ \ t>0.
\end{equation}
\end{prop}

\proof\\
The second inequality in \eq{6.10} readily follows from the comparison
principle and the properties of the linear heat equation since $u$ is
a subsolution to the linear heat equation. As for the first
inequality, we proceed as follows: consider $n\ge 1$ large enough such
that 
\begin{equation}
\label{bof1}
n\ge K_1^{-1} \;\;\mbox{ and }\;\; n^{p-1} > 4 \vert a\vert C_0^2
I_1(p) , 
\end{equation}
the parameters $C_0$, $I_1(p)$ and $K_1$ being defined in \eq{4.6},
\eq{6.1} and Proposition~\ref{fantasio}, respectively. By \eq{bof1},
we have
$$
\frac{1}{2 C_0 n} > \frac{2 \vert a\vert C_0 I_1(p)}{n^p}
$$
and we fix
\begin{equation}
\label{bof2}
\beta_n \in \left( \frac{2 \vert a\vert C_0 I_1(p)}{n^p} , \frac{1}{2
C_0 n} \right).
\end{equation}

We next take the initial datum in \eq{part4.1.1} to be
$u_{0,n}=\beta_n e_1$ and denote by $u_n$ the corresponding solution
to \eq{part4.1.1}. Owing to \eq{1.3} and \eq{4.6} we have
\begin{eqnarray*}
u_n(t) & \ge & \beta_n e^{t\Delta}e_1 - \vert a\vert \int\limits_0^t
\Vert e^{(t-s)\Delta}|\nabla u_n(s)|^p \Vert_\infty \,ds\\
& \geq & \beta_n e^{-t\lambda_1}e_1 -|a| C_0 \int\limits_0^t
e^{-\lambda_1(t-s)} \|\nabla u_n(s)\|^p_{\infty}\,ds.
\end{eqnarray*}
Since $\Vert u_{0,n}\Vert_\infty=\beta_n\le 1/(2C_0 n)$ and $1/n\le
K_1$, we deduce from Proposition~\ref{fantasio} that $u\in
X_{1/n}$. Consequently,
$$\di u_n(t)\geq \beta_n e^{-t\lambda_1}e_1 - |a| C_0 n^{-p}
\int\limits_0^te^{-\lambda_1(t-s)}(1+s^{-1/2})^p e^{-\lambda_1 ps}\,
ds \ge e^{-t\lambda_1} \left( \beta_n e_1 - \frac{|a| C_0 I_1(p)}{n^p}
\right),$$
whence, since $\Vert e_1\Vert_\infty=1$ and $\beta_n$ fulfils \eq{bof2},
\begin{equation}\label{6.13}
\|u_n(t)\|_{\infty}\geq e^{-t \lambda_1} \left( \beta_n-\frac{|a| C_0
I_1(p)}{n^p} \right) \ge e^{-t \lambda_1} \frac{\beta_n}{2}= e^{-t
\lambda_1} \frac{\Vert u_{0,n}\Vert_\infty}{2}.
\end{equation}

Consider now $u_0$ as in Proposition~\pr{prop.1} and denote by $u$ the
corresponding classical solution to \eq{part4.1.1}. Since $\beta_n\to 0
$ as $n\to\infty$, there is $n_0$ fulfilling \eq{bof1} such that
$u_0\ge \alpha e_1 \ge \beta_{n_0} e_1$. By the comparison principle,
we have $u(t)\ge u_{n_0}(t)$ and Proposition~\ref{prop.1} is then a
straightforward consequence of \eq{6.13}.
\endproof

%%%%%%%%%%%%%%%%%%%%%%%%%%%%%%%%%%%%%%%%%%%%%%%%%%%%%%%%%%%%%%%%%%%%%%%%%%

\section{Proof of Theorem \pr{t.7}}
We first consider the case $a<0$. On the one hand, $u$ is a
subsolution of the linear heat equation with the same initial
datum and the comparison principle entails that
\begin{equation}
\label{z1}
0\leq u\leq v \;\;\mbox{ in }\;\; Q_{\infty} \;\;\mbox{ with }\;\;
v(t)=e^{t\Delta}u_0\,,\quad t\ge 0\,.
\end{equation}
On the other hand, by the maximum principle
(see \cite{F,Li} and also \cite[Remark~3.3]{S}) we have
\begin{equation}\label{4.15}
\|\nabla u(t)\|_{\infty} = \|\nabla u(t)\|_{\partial\Omega,\infty}
= \|u_\nu(t)\|_{\partial\Omega,\infty},
\end{equation}
and
\begin{equation}\label{4.15s}
\|u_0\|_{C^1(\overline{\Omega})} \geq \|\nabla v(t)\|_{\infty}
= \|\nabla v(t)\|_{\partial\Omega,\infty}
= \|v_\nu(t)\|_{\partial\Omega,\infty},
\end{equation}
the last equality in \eq{4.15} and \eq{4.15s} being a consequence of
the homogeneous Dirichlet boundary conditions. Owing to the homogeneous
Dirichlet boundary conditions, we deduce from \eq{z1} that
$v_\nu(t)\leq u_\nu(t)\leq 0$ for $t\ge 0$, which, together with
\eq{4.15} and \eq{4.15s}, yields
\begin{equation}\label{4.16}
\|\nabla u(t)\|_\infty = \|u_\nu(t)\|_{\partial\Omega,\infty} \leq
\|v_\nu(t)\|_{\partial\Omega,\infty}=\|\nabla
v(t)\|_{\infty} \leq \|u_0\|_{C^1(\overline{\Omega})}\ \ \
\mbox{for}\ \ \ t>0.
\end{equation}
Consequently, since $p>2$, 
$$u_t-\Delta u=-|a||\nabla u|^{p-2}|\nabla u|^2\geq
-|a|\|u_0\|_{C^1(\overline{\Omega})}^{p-2}|\nabla u|^2 \ \ \
\mbox{in} \ \ \ Q_{\infty}.$$
Denote by $C_6=|a|\|u_0\|_{C^1(\overline{\Omega})}^{p-2}$ and let $w$
be the solution to the following initial boundary value problem:
\begin{equation}\label{4.17}
\left\{\begin{array}{lll} \di w_t-\Delta w=-C_6|\nabla w|^2 \
\
\mbox{in}\ \ \ Q_{\infty},\\
w(t,x)=0 \ \ \mbox{on}\ \  \Gamma_{\infty},\\
w(0,.)=u_0  \ \ \mbox{in} \ \ \Omega.\end{array}\right.
\end{equation}
Then $z=1-e^{-C_6w}$ satisfies the linear heat equation in
$Q_{\infty}$ with homogeneous Dirichlet boundary conditions and
initial datum  $z(0)=1-e^{-C_6 u_0}$, so that
$w(t)=-C_6^{-1}\log\left(1-e^{t\Delta}z(0)\right)$ for $t\ge 0$,
while the comparison principle ensures that
$$-\frac{1}{C_6}\log\left(1-e^{t\Delta}z(0)\right)\leq
u(t)\leq e^{t\Delta} u_0, \ \ \ \mbox{for}\ \ \ t>0.$$
Since $-\log(1-r)\ge r$ for $r\in (0,1)$ and $e^{t\Delta}z(0)\in (0,1)$
for $t\ge 0$, we conclude that
$$ \frac{1}{C_6} e^{t\Delta}z(0)\leq u(t)\leq e^{t\Delta} u_0,
\ \ \ \mbox{for}\ \ \ t>0,$$
and \eq{4.9} holds true with $w_0=(1-e^{-C_6 u_0})/C_6$
which is a positive function in $C_0^1(\overline{\Omega})$ and $W_0=u_0$.
Furthermore, the large time behaviour of $\nabla u(t)$ is a
consequence of \eq{4.7} and \eq{4.16}:
$$\|\nabla u(t)\|_{\infty}\leq\|\nabla e^{t\Delta} u_0\|_{\infty}\leq
C_0(1+t^{-1/2})e^{-\lambda_1 t}\|u_0\|_{\infty},$$
and the proof of Theorem~\pr{t.7} is complete for $a<0$.

We now turn to the case $a>0$. By \cite[Proposition 3.1]{S}, the
condition $\|u_0\|_{C^1(\overline{\Omega})} \leq \eps$ not only
warrants that the corresponding classical solution $u$ to
\eq{part4.1.1} is global but also that it is bounded in
$C^1(\overline{\Omega})$. Consequently, there is a positive constant
$C_7>0$ such that
$$\|\nabla u(t)\|_{\infty}\leq C_7 \ \ \mbox{for}\ \ t\geq
0.$$
Thanks to this property, we may proceed as in the previous case and
deduce from the comparison principle that
\begin{equation}\label{5.1}
e^{t\Delta} u_0 \leq u(t) \leq w(t) = \frac{1}{C_8}
\log(1+e^{t\Delta}z(0)) \ \ \mbox{for}\ \ t>0,
\end{equation}
where $C_8=a C_7^{p-2}$, $z(0)=e^{C_8 u_0}-1$ and $w$ is the solution
to the following initial boundary value problem:
\begin{equation}\label{4.18}
\left\{\begin{array}{lll} \di w_t-\Delta w=C_8|\nabla w|^2 \ \
\mbox{in}\ \ \ Q_{\infty},\\
w(t,x)=0 \ \ \mbox{on}\ \  \Gamma_{\infty},\\
w(0,.)=u_0  \ \ \mbox{in} \ \ \Omega.\end{array}\right.
\end{equation}
Since $\log(1+r)\leq r$ for $r\ge 0$ and $e^{t\Delta} z(0)\ge 0$ for
$t\ge 0$, we infer from \eq{5.1} that \eq{4.9} is satisfied with
$w_0=u_0$ and $W_0=(e^{C_8 u_0}-1)/C_8$. In addition, owing to
\eq{5.1} and \cite[Remark 3.3]{S} we have that:
\begin{eqnarray*}
& & \|\nabla u(t)\|_{\infty} = \|\nabla u(t)\|_{\partial\Omega,\infty}
= \|u_\nu(t)\|_{\partial\Omega,\infty} \\
& \leq & \|w_\nu(t)\|_{\partial\Omega,\infty} = \|\nabla
w(t)\|_{\partial\Omega,\infty} = \|\nabla w(t)\|_{\infty} =
\frac{1}{C_8} \left\| \frac{\nabla e^{t\Delta} z(0)}{1+
e^{t\Delta}z(0)} \right\|_\infty,
\end{eqnarray*}
from which the estimate \eq{4.5} follows.
\endproof

%%%%%%%%%%%%%%%%%%%%%%%%%%%%%%%%%%%%%%%%%%%%%%%%%%%%%%%%%%%%%%%%%%%%%%

\section{Proof of Theorem~\pr{t.9}}
We first prove Theorem~\pr{t.9} for non-negative initial data $u_0\in
C_0([-1,1])$ which are profiled, that is, $u_0$ is a non-decreasing
function on $(-1,0)$ and a non-increasing function on $(0,1)$. From
\cite[Corollary 4.4]{GGK} we know that this property is preserved
throughout time evolution, so that $u(t)$ is a non-decreasing function
on $(-1,0)$ and a non-increasing function on $(0,1)$ for any $t>0$.
Since $u(t)\in C^1([-1,1])$ for $t>0$, an alternative formulation of
this property is $|u_x(t,x)|=- \mbox{sign}(x) u_x(t,x)$ for $(t,x)\in
(0,\infty)\times (-1,1)$. Therefore, $u$ also solves
\begin{equation}\label{7.5}
\left\{\begin{array}{lll} u_t-u_{xx}=-a u_x & \mbox{in} &
(0,\infty)\times(0,1),\\
u_x(t,0)=u(t,1)=0 & \mbox{in} &(0,\infty),\\
u(0,x) = u_0(x)& \mbox{in} &(0,1),\end{array}\right.
\end{equation}
and a similar equation on $(-1,0)$. Conversely, as a consequence of
the uniqueness of the solution to \eq{part4.1.1}, solving \eq{7.5}
on $(0,1)$ and $(-1,0)$ gives back the solution to \eq{part4.1.1}.
We shall therefore study the solution to \eq{7.5}. Using the
transformation
\begin{equation}\label{7.6}
v(t,x)=e^{a^2t/4}e^{-ax/2}u(t,x)\,, \quad (t,x)\in
(0,\infty)\times (0,1)\,,
\end{equation}
then $v$ satisfies the following problem
\begin{equation}\label{7.7}
\left\{\begin{array}{lll} v_t-v_{xx}=0& \mbox{in} &
(0,\infty)\times(0,1),\\
2v_x(t,0)+av(t,0)=v(t,1)=0 & \mbox{in} & (0,\infty),\\
v(0,x)=v_0(x)=e^{-ax/2}u_0(x) & \mbox{in} & (0,1),
\end{array}\right.
\end{equation}
which also reads $v_t = L v$ with $v(0)=v_0$, the unbounded linear
operator $L$ being defined in Theorem~\pr{t.9}. The initial boundary
value problem \eq{7.7} being linear, the large time behaviour of its
solutions is determined by the spectrum of $L$. First, classical results
ensure that the spectrum is an increasing sequence $(\alpha_n)_{n\ge 1}$
of eigenvalues converging to $\infty$ and the corresponding normalized
eigenfunctions $(\varphi_n)_{n\ge 1}$ form an orthonormal basis of
$L^2(0,1)$ (see, e.g., \cite[Th\'eor\`eme~6.2-1 and Remarque~6.2-2]{RT}).
The next step is to identify the eigenvalues and eigenfunctions of $L$.

\begin{prop}\label{p.20}
For $a\ne 0$, the equation $\tan(z)=2z/a$ has a countably infinite
number of positive solutions and we denote by $\mathcal{Z}_a$ the set
of these solutions.

\begin{itemize}

\item[$(i)$] If $a\in(-\infty,2)\setminus\{0\}$, we have $\{ \sqrt{\alpha_n}
; n\ge 1\} = \mathcal{Z}_a$ with $\sqrt{\alpha_n}\in \left(
((2n-1)\pi)/2,n\pi \right)$ if $a<0$ and $\sqrt{\alpha_n}\in \left(
(n-1)\pi,((2n-1)\pi)/2 \right)$ if $a\in (0,2)$ for $n\geq 1$.

\item[$(ii)$] If $a\in[2,\infty)$, we have $\{ \sqrt{\alpha_n} ; n\ge 2\} =
\mathcal{Z}_a$ and $\sqrt{\alpha_n}\in \left( (n-1)\pi,((2n-1)\pi)/2
\right)$ for $n\geq 2$. In addition, if $a=2$ then $\alpha_1=0$ and
the corresponding eigenfunction is given by $\varphi_1(x) =
\sqrt{3}(1-x).$

\noindent If  $a>2$ then $(-\alpha_1)$ is the unique positive real
number satisfying
$$
\frac{a+2\sqrt{-\alpha_1}}{a-2\sqrt{-\alpha_1}}=e^{2\sqrt{-\alpha_1}}
\;\;\mbox{ with }\;\; \alpha_1\in \left(
-\frac{a^2}{4},-\frac{a(a-2)}{4} \right), 
$$
and the corresponding eigenfunction is given by $\varphi_1(x) = A_1(a)
\sinh{\left( \sqrt{-\alpha_1}(1-x) \right)}$, the parameter $A_1(a)$
being a positive constant such that $\|\varphi_1\|_{2}=1$ and
$\varphi_1>0$ in $(0,1)$.

\end{itemize}

Moreover, for $n\ge 1$ such that $\sqrt{\alpha_n}\in\mathcal{Z}_a$,
the corresponding eigenfunction $\varphi_n$ is given by
$\varphi_n(x)=A_n(a) \sin(\sqrt{\alpha_n}(1-x))$ where $A_n(a)$ is
chosen such that $\|\varphi_n\|_{2}=1$. If $n=1$, we also choose
$A_1(a)$ such that $\varphi_1>0$ in $(0,1)$.
\end{prop}

\begin{rem}\label{r.21}
By Proposition~\ref{p.20}, all the eigenvalues of $L$ are positive if
$a<2$, $a\ne 0$. A direct proof of this fact can be performed as
follows: let $\alpha$ be an eigenvalue of $L$ with corresponding
eigenfunction $\varphi$, so that $L\varphi=\alpha\varphi$. Multiplying
this identity by $\varphi$ and integrating over $(0,1)$ we have:
$$\int\limits_0^1 \vert\varphi_x(x)\vert^2dx -
\frac{a}{2}\varphi(0)^2=\alpha\int\limits_0^1\varphi(x)^2\, dx.$$
Since $\varphi(1)=0$, an elementary computation shows that
$$\varphi(0)^2=\left(\int\limits_0^1\varphi_x(x)dx\right)^2 \leq
\int\limits_0^1\vert\varphi_x(x)\vert^2dx,$$
so that the left-hand side of the above identity is positive if $a<2$.
\end{rem}

{\it Proof of Proposition~\pr{p.20}:}\\
Let $\alpha$ be an eigenvalue of $L$ with corresponding eigenfunction
$\varphi$. Then
\begin{equation}
\label{7.3}
-\varphi_{xx} = \alpha\varphi \;\;\mbox{ in }\;\; (0,1) \;\;\mbox{ and
}\;\; 2\varphi_x(0)+a\varphi(0)=\varphi(1)=0\,.
\end{equation}

In order to solve \eq{7.3} we distinguish among the cases $\alpha<0$,
$\alpha=0$ and $\alpha>0$.

$1)$ If $\alpha<0$, solving the first equation in \eq{7.3} gives
$$
%\begin{equation}\label{10.2}
\varphi(x)=Ae^{x\sqrt{-\alpha}}+Be^{-x\sqrt{-\alpha}},\ \ \
x\in(0,1),
$$
%\end{equation}
for some yet unspecified real numbers $A$ and $B$. To comply with the
boundary conditions in \eq{7.3}, we deduce that $\alpha$ has to verify
the following equation:
\begin{equation}\label{10.1}
\frac{a+2\sqrt{-\alpha}}{a-2\sqrt{-\alpha}}=e^{2\sqrt{-\alpha}}.
\end{equation}
Now, it is easy to check that the equation $e^z =(a+z)/(a-z)$ has a
unique positive solution $\varrho_a$ if and only if $a>2$, and
$\varrho_a\in \left( \sqrt{a(a-2)},a \right)$. Consequently, if $a>2$,
we have $\alpha_1=-\varrho_a^2/4\in (-a^2/4,-a(a-2)/4)$ and
$\varphi_1(x) = A_1(a) \sinh\left( \sqrt{-\alpha_1}(1-x) \right)$ for
$x\in(0,1)$ with
\begin{equation}\label{10.11}
A_1(a) = 2\ \left( \frac{\sinh\left( 2\sqrt{-\alpha_1}
\right)}{\sqrt{-\alpha_1}} - 2 \right)^{-1/2},
\end{equation}
so that
$\|\varphi_1\|_{2}=1$ and $\varphi_1$ is positive in $(0,1)$.

$2)$ For \eq{7.3} to have a non-zero solution with $\alpha=0$, it is
necessary that $a=2$. Hence, if $a=2$, we have $\alpha_1=0$ with
$\varphi_1(x)=\sqrt{3}(1-x)$, $x\in (0,1)$.

$3)$ If $\alpha>0$, solving the first equation in \eq{7.3} leads to
$$
%\begin{equation}\label{10.3}
\varphi(x)=A\sin(\sqrt{\alpha}x)+B\cos(\sqrt{\alpha}x),\ \ \
x\in(0,1),
$$
%\end{equation}
for some yet unspecified real numbers $A$ and $B$. Requiring that
$\varphi$ fulfils the boundary conditions in \eq{7.3} implies that
$\alpha\in\mathcal{Z}_a$. Then either $a\ge 2$ and, since $\alpha_1$
has already been determined, we have $\{ \sqrt{\alpha_n} ; n\ge 2\}
= \mathcal{Z}_a$. Or $a<2$ ($a\ne 0$) and $\{ \sqrt{\alpha_n} ;
n\ge 1\} = \mathcal{Z}_a$. In both cases, $\varphi_n(x) = A_n(a)
\sin(\sqrt{\alpha_n}(1-x))$ for $x\in(0,1)$ with
\begin{equation}\label{10.10}
A_n(a) = 2\ \left( 2-\frac{\sin\left( 2\sqrt{\alpha_n}
\right)}{\sqrt{\alpha_n}} \right)^{-1/2},
\end{equation}
chosen such that $\|\varphi_n\|_{2}=1$ and $\varphi_1$ is positive in
$(0,1)$.
\endproof

\medskip

As a direct consequence of formulae \eq{10.11} and \eq{10.10}, we
next derive some properties of $(A_n(a))_{n\ge 1}$ according to the
values of $a$.

\begin{lemma}
\label{coeff}
If $a\ne 0$ and $n\ge 2$, we have
\begin{equation}\label{10.12}
A_n(a) \leq \sqrt{\pi},
\end{equation}
and 
$$
\lim_{a\to 0} A_1(a)=\sqrt{2}\,, \quad \lim_{a\to 2} A_1(a)=\infty\,,
\quad \lim_{a\to\infty} A_1(a)=0\,.
$$
Moreover, if $a<0$, $A_1(a)\le\sqrt{\pi}$.
\end{lemma}

{\it Proof of Theorem~\ref{t.9}:}\\
Since the normalised eigenfunctions $(\varphi_n)_{n\geq 1}$ of $L$ form
an orthonormal basis in $L^2(0,1)$, the solution to \eq{7.7} is given by
$$
v(t,x)=\sum\limits_{n\geq 1} e^{-\alpha_n t}
<v_0,\varphi_n>_{L^2(0,1)} \varphi_n(x) \ \ \mbox{in} \ \ \
(0,\infty)\times(0,1),
$$
where $<.,.>_{L^2(0,1)}$ denotes the usual scalar product in
$L^2(0,1)$. From \eq{7.6} we deduce that
\begin{equation}\label{7.8}
u(t,x)=\sum\limits_{n\geq 1} e^{-((a^2/4)+\alpha_n) t} e^{ax/2}
<v_0,\varphi_n>_{L^2(0,1)} \varphi_n(x) \ \ \mbox{in} \ \ \
(0,\infty)\times(0,1).
\end{equation}
Changing $x$ to $-x$ we obtain a similar identity on the interval
$(-1,0)$
\begin{equation}\label{7.9}
u(t,x)=\sum\limits_{n\geq 1} e^{-((a^2/4)+\alpha_n)t} e^{-ax/2}
<\tilde{v}_0,\varphi_n>_{L^2(0,1)} \varphi_n(-x) \
\ \mbox{in} \ \ \ (0,\infty)\times(-1,0),
\end{equation}
where $\tilde{v}_0(y)=e^{-ay/2}u_0(-y)$ for $y\in(0,1)$.

Thanks to the properties of the eigenvalues $(\alpha_n)_{n\geq 1}$
and to relations \eq{10.11}, \eq{10.10} and \eq{10.12}, a simple
computation shows that, if $a>0$,
\begin{eqnarray*}
\|u(t)\|_{\infty} & \le & e^{a/2} \|v_0\|_{\infty} \sum\limits_{n\ge
1} e^{-((a^2/4)+\alpha_n)t} \|\varphi_n\|_{\infty}^2\\
& \le & e^{a/2} \|u_0\|_{\infty} e^{-((a^2/4)+\alpha_1)t} \left(
\|\varphi_1\|_{\infty}^2+\sum\limits_{n\ge 2}
e^{-(\alpha_n-\alpha_1)t} A_n(a)^2 \right)\\ 
& \le & e^{a/2} \|u_0\|_{\infty} e^{-((a^2/4)+\alpha_1)t} \left(
\|\varphi_1\|_{\infty}^2 + \pi\sum\limits_{n\ge 2}
e^{-(2n-3)(2n-1)\pi^2 t/4} \right)\\ 
& \le & e^{a/2} \|u_0\|_{\infty} e^{-((a^2/4)+\alpha_1)t} \left(
\|\varphi_1\|_{\infty}^2 + \pi\sum\limits_{n\ge 1} e^{-n\pi^2t/2} \right)\\
& \le & e^{a/2} \left( \|\varphi_1\|_{\infty}^2 +
\frac{\pi}{e^{\pi^2t/2}-1} \right) \|u_0\|_{\infty} e^{-((a^2/4)+\alpha_1)t},
\end{eqnarray*}
for $t>0$, whence \eq{7.2} for $t\ge 1$.

If $a<0$ we notice that \eq{10.12} holds for all $n\geq 1$. Therefore,
we have
\begin{eqnarray*}
\|u(t)\|_{\infty} & \le & \|v_0\|_{\infty} \sum\limits_{n\ge 1}
e^{-((a^2/4)+\alpha_n)t} A_n(a)^2\\
& \le & e^{|a|/2} \|u_0\|_{\infty} e^{-((a^2/4)+\alpha_1)t} \pi \left(
1+\sum\limits_{n\ge 2} e^{-(\alpha_n-\alpha_1)t} \right)\\
& \le & e^{|a|/2} \frac{\pi e^{\pi^2t/2}}{e^{\pi^2t/2}-1}
\|u_0\|_{\infty} e^{-((a^2/4)+\alpha_1)t},
\end{eqnarray*}
whence \eq{7.2} for $a<0$. We have thus established \eq{7.2} for a
profiled function $u_0$.

\medskip

In the general case, if $u_0\in C_0^+([-1,1])$, we define
$\bar{u}_0$ by:
$$\bar{u}_0(x)=\sup\{\ u_0(y); \ \ \vert y\vert \ge |x| \ \}.$$
Thus, $\bar{u}_0\in C_0([-1,1])$ is a profiled function such that
$u_0\leq\bar{u_0}$. Denoting by $\bar{u}$ the solution to
\eq{part4.1.1} corresponding to the initial datum $\bar{u}_0$, we
infer from the comparison principle and the estimate \eq{7.2} for
$\bar{u}$ that
$$\|u(t)\|_{\infty}\leq\|\bar{u}(t)\|_{\infty} \leq \gamma(a)\
\|\bar{u}_0\|_{\infty}\ e^{-((a^2/4)+\alpha_1)t}, \ \ \ t\ge 1.$$
Since $\|u_0\|_{\infty}=\|\bar{u}_0\|_{\infty}$, we deduce
that \eq{7.2} is fulfilled in the general case too.
\endproof

The following corollary is a direct consequence of formulae \eq{7.8}
and \eq{7.9}.

\begin{cor}\label{c.1}
Let $u_0\in C_0^+(\overline{\Omega})$ be a profiled function. Under
the hypotheses of Theorem~\pr{t.9}, the solution $u$ to \eq{part4.1.1}
also satisfies
\begin{equation}\label{7.20}
\left| e^{((a^2/4)+\alpha_1)t} e^{-ax/2} u(t,x) -
<v_0,\varphi_1>_{L^2(0,1)} \varphi_1(x) \right| \le \gamma(a)\
\|u_0\|_{\infty}\ e^{-(\alpha_2-\alpha_1)t}
\end{equation}
for $t\ge 1$ and $x\in (0,1)$ and a similar inequality on $(-1,0)$ with
$\tilde{v}_0$ instead of $v_0$.
\end{cor}

As a final comment, we emphasize that the large time behaviour of
solutions to \eq{part4.1.1} is rather peculiar in the case $p=1$ and
$N=1$, since it is the only situation where we observe a real
difference between the solution to the linear heat equation and the
solution to \eq{part4.1.1}. The nonlinear term actually plays an
important role whatever the sign of $a$ is. Indeed, recalling that the
first eigenvalue $\lambda_1$ of the Laplace operator with homogeneous
Dirichlet boundary conditions is given by $\lambda_1=\pi^2/4$ in the
particular case $\Omega=(-1,1)$, we denote by
$r_1(a)=(a^2/4)+\alpha_1$ the exponent which gives the decay rate in
\eq{7.2}, $\alpha_1$ being the first eigenvalue of the operator $L$
defined in Theorem~\pr{t.9} and thus depending on $a$. We then aim at
comparing $r_1(a)$ and $\lambda_1$. 

$(i)$ if $a<0$ we have $r_1(a)>\lambda_1$ and the absorption term
$-\vert a\vert \vert u_x\vert$ drives the solutions to \eq{part4.1.1}
to zero at a faster rate than the solutions to the linear heat
equation. Moreover, we have $r_1(a)\searrow\lambda_1$ as $a\nearrow 0$
and $r_1(a)\nearrow \infty$ as $a\searrow -\infty$. 

$(ii)$ if $a>0$, we have $r_1(a)\in (0,\lambda_1)$ and the source term
$a \vert u_x\vert$ slows down the convergence to zero of solutions to
\eq{part4.1.1}. Furthermore, $r_1(a)\nearrow\lambda_1$ as $a\searrow 0$ and
$r_1(a)\searrow 0$ as $a\nearrow\infty$. Indeed, if $a\in (0,2)$, we
have $\sqrt{\alpha_1}\in (0,\pi/2)$ and
$tan{(\sqrt{\alpha_1})}=2\sqrt{\alpha_1}/a$ by
Proposition~\pr{p.20}~(i), whence
$$
r_1(a)=\left( \frac{\sqrt{\alpha_1}}{\sin(\sqrt{\alpha_1})}\right)^2
\le \frac{\pi^2}{4}=\lambda_1,
$$
while, for $a>2$, it follows from Proposition~\pr{p.20}~(ii) and
\eq{10.1} that $\alpha_1\in (-a^2/4,-a(a-2)/4)$ and 
$$
r_1(a)=-\frac{4\alpha_1e^{2\sqrt{-\alpha_1}}}{(e^{2\sqrt{-\alpha_1}}-1)^2}
\leq 1 \le \lambda_1.
$$
Finally, by \eq{10.1}, $r_1(a)$ also satisfies
$$
r_1(a)=\frac{1}{4}(a+2\sqrt{-\alpha_1})(a-2\sqrt{-\alpha_1})=
\frac{1}{4}(a+2\sqrt{-\alpha_1})^2e^{-2\sqrt{-\alpha_1}},$$
from which we deduce that $r_1(a)\searrow 0$ as $a\nearrow\infty$.
\bigskip

{\bf Acknowledgements.} We would like to thank Professor Philippe
Souplet for interesting discussions related to this work and Professor
Guy Barles for suggesting to us a shorter proof of
Proposition~\ref{part4.p.1}.

%\bibliographystyle{plain}
%\bibliography{bibliografie}

\end{document}